\documentclass{academic}
\usepackage{amssymb,amsfonts,amsbsy}
\newtheorem{example}[theorem]{Example}
\newtheorem{con}[theorem]{Convention}

\def\KK{{\mathbb K}} 
 
\def\kk{{\overline{k}}}

\def\II{{\mathbb I}}
\def\ii{{\mathcal I}}
 
\def\LL{{\mathcal L}} 
\def\MM{{\mathcal M}} 
\def\OO{{\mathcal O}}

\def\R{{\mathcal R}}
\def\CC{{\mathbb C}}
\def\TT{{\mathbb T}}
\def\PP{{\mathbb P}}

\def\QQ{{\mathbb Q}}

\def\ZZ{{\mathbb Z}}
\newcommand{\Rd}{{\rm Res}}

\begin{document}
\shortauthor{C. D'Andrea
}
\shorttitle{Resultants and Moving Surfaces}
\title{Resultants and Moving Surfaces}
\author{Carlos D'Andrea\footnote{\noindent Partially supported by Universidad de Buenos Aires, grant TX094, and
Agencia Nacional de Promoci\'on Cient\'{\i}fica y Tecnol\'ogica
(Argentina), grant 3-6568.}}
\address{ Departamento  de Matem\'atica, F.C.E.y N.,  
Universidad de Buenos Aires. Pabell\'on I, Ciudad Universitaria 
(1428) Buenos Aires, Argentina}
\maketitle

\abstract{We prove a conjectured relationship among  resultants and the determinants 
arising in the formulation of the method  of moving surfaces for computing
the implicit equation of rational surfaces formulated by Sederberg.  
In addition, we extend the validity of this method
to the case of not properly parametrized surfaces without base
points.}

\section{Introduction}
Given four polynomials in two variables $x_1(s,t), \ x_2(s,t),$ \
$x_3(s,t)$ and $ x_4(s,t),$ the equations
\begin{equation}
\label{parametric}
X_1=\frac{x_1(s,t)}{x_4(s,t)}, \ \ X_2=\frac{x_2(s,t)}{x_4(s,t)}, \ \
X_3=\frac{x_3(s,t)}{x_4(s,t)}, 
\end{equation}
define a parametrization of a rational surface.
\par
The implicitization problem consists in finding  another polynomial 
$F(X_1,X_2,X_3)$ such that 
$F(X)=0$ is the equation of the smallest algebraic surface 
containing (\ref{parametric}).
\par
A classical method for finding this implicit equation is to eliminate the
variables $s$ and $t$ by computing 
the bivariate resultant of the polynomials
\begin{equation}
\label{ppert}
x_1(s,t)-X_1\,x_4(s,t),\ \  x_2(s,t)-X_2\,x_4(s,t), \ \
x_3(s,t)-X_3\,x_4(s,t).
\end{equation}
There are several types of bivariate resultant \cite{dix,GKZ,stu}. They are related with different  compactifications of the
affine space, where  the input polynomials are defined.
For example, if we view the polynomials $x_i$ as general polynomials of degree less 
than or  equal to 
$n,$ the {\it multivariate resultant} 
may be taken, and its vanishing means that 
the system (\ref{ppert}) has a common root in the projective space $\PP^2.$
This is the so called {\it dense} or {\it triangular case.}
\par
Another situation is when we regard the polynomials $x_i$ as having degree less than or 
equal to $m$ in $s,$ and less than or equal to $n$ in $t.$  Here, one can use 
the {\it bihomogeneous resultant.} This is the {\it 
tensor product case,} and the vanishing of the resultant means
that the system has a  solution in $\PP^1\times\PP^1.$
\par
In both cases, the resultant of (\ref{ppert}) gives the implicit equation 
(actually, a power of it) in 
the absence of {\it base
points} i.e. when there are no $(s_0, t_0)$ in the corresponding space 
($\PP^2$ or $\PP^1\times\PP^1$) such that
$$x_1(s_0,t_0)=x_2(s_0,t_0)=x_3(s_0,t_0)=x_4(s_0,t_0)=0.$$
In \cite{SC}, Sederberg and Chen introduced a new technique called
{\it the moving quadric method} for finding the
implicit equation of (\ref{parametric}). It uses smaller determinants than
the classical methods and often works in the presence of base points.
\par
A detailed analysis of this technique is given in \cite{CGZ}, where 
sufficient conditions are established for the validity of 
implicitization
by the method of moving quadrics for rectangular tensor product surfaces and 
triangular surfaces in the absence of base points, and when the surface
is properly parametrized, i.e. in when the parametrization is not necessarily
one-to-one.
In that paper the authors conjectured a relationship between the moving
plane and the moving quadric coefficient matrices for both the tensor product
and the triangular case (Conjectures $6.1$ and $6.2$ in \cite{CGZ}). 
The intuition behind this conjecture was a similar relationship valid in
the plane case ( see for instance, \cite{ZCG}).
\par
This paper presents a general relationship between the moving plane and
the moving surfaces of degree $d$ coefficient matrices. As a special case, 
when $d=2,$ it provides a proof for both conjectures. In addition, the validity
of the method when the surface has not base points, but is not 
necessarily properly parametrized, is proven.
\par
The approach consists in the factorization of the moving surface coefficient 
matrices as composition of simpler matrices of linear maps. 
This, combined with the well-known 
formulation of
the resultant as the determinant of a Koszul complex \cite{GKZ}, gives the
desired results which recover, as a particular case, Theorems $4.1$ and
$5.1$ in \cite{CGZ}. Adapting these techniques to the planar case,
one can also produce alternative proofs of similar relationships 
given in \cite{SGD} and \cite{ZCG}.
\par
The paper is organized as follows:  in Section \ref{predos}, some geometric 
definitions are established, in order
to provide a better understanding of the relationships  which will follow. In 
Section \ref{dos},  the relations between 
maximal minors of moving surfaces matrices and the resultant are proven for 
surfaces parametrized by bihomogeneous polynomials. 
In the following section, the validity of the method of moving quadrics is
extended to the case of surfaces without base points but not properly 
parametrized. 
Finally, in the last section, the same situation is considered for 
surfaces parametrized by homogeneous polynomials.
\par
In closing, it should be mentioned that there are recent results concerning the implicitization
problem in the presence of basepoints which have been obtained via the study of syzygies in \cite{C}, and
the use of residual resultants in \cite{B}.
\bigskip
\section{Moving Surfaces}\label{predos}
This section will  review some  basic notions used in the ``method 
of moving conics and quadrics'', as stated in 
\cite{CGZ,SC,ZCG}, in order to  provide a geometric meaning of the algebraic 
tools to be developed in the 
following paragraphs.
\par
Let $\KK$ be a field.
A $d$-surface is an implicit homogeneous equation in the variables 
$X_1,X_2,X_3$ and $X_4$ of degree $d:$
$$
\sum_{|\gamma|=d} c_\gamma X^\gamma=0, \ c_\gamma\in\KK.
$$
A moving $d$-surface of bi-degree $(\sigma_1,\sigma_2)$ is a family of
$d$-surfaces parametrized by $s,u,t$ and $v$ as follows:
\begin{equation}
\label{msurface}
\sum_{i=0}^{\sigma_1}\sum_{j=0}^{\sigma_2}{
\left(\sum_{|\gamma|=d} A^{ij}_{\gamma} X^\gamma\right)s^i\,u^{\sigma_1-i}\,
t^{j}\,v^{\sigma_2-j}}=0, \ A^{ij}_\gamma\in\KK.
\end{equation}
For each fixed value of the parameters, equation (\ref{msurface}) is an
implicit equation of degree $d$ in $\KK^3.$
\par
Similarly, a moving $d$-surface of degree $\sigma,$ is defined by:
\begin{equation}
\label{msurface2}
\sum_{i+j\leq \sigma}{
\left(\sum_{|\gamma|=d} A^{ij}_{\gamma} X^\gamma\right)s^i\,
t^{j}\,u^{\sigma-i-j}}=0, \ A^{ij}_\gamma\in\KK.
\end{equation}
In both cases, a moving $1$-surface will be called ``moving
plane''. If $d=2,$ it is a  ``moving quadric'' (cf. \cite{CGZ}).
\par
Given a family of four bihomogeneous (resp. homogeneous) polynomials
$ x_i (s,u;v,t)$ (resp. $ x_i(s,t,u)$) of bi-degree $(m,n)$ (resp. degree $n$) 
with coefficients in $\KK,$ the moving $d$-surface
(\ref{msurface}) (resp. (\ref{msurface2})) is said to follow
the rational surface $\left(\frac{ x_1}{ x_4},
\frac{ x_2}{ x_4}, \frac{ x_3}{ x_4}
\right)$ if 
$$
\sum_{i=0}^{\sigma_1}\sum_{j=0}^{\sigma_2}{
\left(\sum_{|\gamma|=d} A^{ij}_{\gamma}  x^\gamma\right)s^i\,u^{\sigma_1-i}\,
t^{j}\,v^{\sigma_2-j}}=0
$$
resp.
$$\sum_{i+j\leq \sigma}{
\left(\sum_{|\gamma|=d} A^{ij}_{\gamma}  x^\gamma\right)s^i\,
t^{j}\,u^{\sigma-i-j}}=0.
$$
In order to find the $\KK$-vector space of all 
$d$-surfaces that follow the rational surface 
of a fixed bi-degree (resp. degree), 
set the coefficients of all the monomials $s^{\alpha}\,u^{\sigma_1
-\alpha}\,t^\beta\,v^{\sigma_2-\beta}$ (resp. $s^\alpha\,t^\beta\,
u^{\sigma-\alpha-\beta}$) in the implicit equation equal to zero, and solve the
linear system of equations in the indeterminates 
$\{A^{ij}_\gamma\}.$ 
\begin{example}
Consider the following family of homogeneous polynomials:
\begin{equation}
\label{lpm}
\begin{array}{lcl}
{x}_1&=&s^3\\
{x}_2&=&t^3,\\
{x}_3&=&u^3,\\
{x}_4&=&s^3+t^3+u^3.
\end{array}
\end{equation}
They define a parametric surface  contained in the hyperplane
$$X_1+X_2+X_3-X_4=0$$
which is a moving plane of degree $0.$ Note that there exists a
moving plane of degree zero if and only if the surface is contained in
a plane. In this case, this is the only plane which contains (\ref{lpm}), so it
is a basis of the moving planes of degree $0.$ 
Also, it is straightforward to compute a 
family of generators 
for the moving quadrics of the same degree:
$$\begin{array}{ccc}
X_1\left(X_1+X_2+X_3-X_4\right)&=&0,\\
X_2\left(X_1+X_2+X_3-X_4\right)&=&0,\\
X_3\left(X_1+X_2+X_3-X_4\right)&=&0,\\
X_4\left(X_1+X_2+X_3-X_4\right)&=&0.\\
\end{array}$$
\end{example}
\medskip
\begin{example}
\item This example appears in \cite{CGZ}. Set
$$\begin{array}{lcl}
{x}_1&=&s\,t+u\,v,\\
{x}_2&=&s\,v,\\
{x}_3&=&u\,t,\\
{x}_4&=&s\,v+u\,t+u\,v,\\
\end{array}$$
and $\sigma_1=\sigma_2=1.$ A basis of the space of moving planes of
bidegree $(1,1)$ which
follow the parametric surface is given by:
$$\begin{array}{rcl}
(X_4-X_1-X_2-X_3)\,u\,v+s\,v\,X_3&=&0,\\
(X_4-X_1-2\,X_2-X_3)\,u\,v+s\,v(X_4-X_2)&=&0.
\end{array}$$

With the aid of Maple,  a basis of  $24$  moving quadrics of  the same bidegree was
found, $8$ of which
come from the moving planes by multiplication by $X_1,\cdots,X_4.$
\end{example}

\section{ The tensor product surface case}  \label{dos}
\subsection{Notation}
\label{notation}
Let $x_1,\ x_2,\ x_3$ and $x_4$ be four generic bihomogeneous polynomials 
in two variables of bi-degree $(m,n),$ i.e.
$$ x_i\left(s,u;t,v\right) = \sum_{j=0}^m\sum_{k=0}^n
c^i_{jk} \ s^j\,u^{m-j}\,t^k\,v^{n-k} \ \ \ \ i=1,\ldots,4.$$
Set $\KK:= \QQ(c^i_{jk}),$
and let $S_{k,l}$ denote the space of polynomials of bi-degree 
$(k,l)$ with coefficients in $\KK.$
\begin{con}
All spaces  to be considered have a monomial basis. 
Suppose all these bases have a fixed order. 
Then,  matrices ``in the monomial bases'' may be defined with
no ambiguity.
\end{con}
\par
Let $\phi$ be the $\KK$-linear map:
\begin{equation}
\label{phi}
\begin{array}{cccc}
\phi:& \mathop{S_{m-1,n-1}}^4& 
\to & S_{2m-1,2n-1}\\
&\left(p_1,p_2,p_3,p_4\right)&\mapsto & \sum_{i=1}^4 p_i\, x_i,\\
\end{array}
\end{equation}
and, following \cite{CGZ}, denote by $MP$ the matrix of $\phi$ in the
monomial bases. It is square, of size $4mn.$

\begin{remark} \
With the definitions stated in the previous section, it is not hard to check that 
$MP$ is the coefficient matrix of the 
linear system generated by
the moving planes of bi-degree $(m-1,n-1)$ that follow the rational
surface given by $\left(\frac{x_1}{x_4},\frac{x_2}{x_4},\frac{x_3}{x_4}\right).$
\end{remark}

Let $d$ be a positive integer and set  
$$\Gamma:=\{\gamma\in{\ZZ_{\geq0}}^4:\,|\gamma| = d\}.$$
Consider the map
$$\Psi^d:\, {S_{m-1,n-1}}^{\Gamma}\, 
\to \, S_{(d+1)m-1,(d+1)n-1}$$ 
which sends
the sequence $\left(p_\gamma\right)_{\gamma\in\Gamma}$ to the polynomial
\begin{equation} 
\label{psi}
\sum_{\gamma\in\Gamma}p_\gamma x^\gamma.
\end{equation}
Let $MQ^d$ be the matrix of $\Psi^d$ in the monomial bases.
Also, set $$\Gamma_0:=\{\gamma\in{\ZZ_{\geq0}}^4:\
|\gamma| = d, \ \gamma_4\leq 1\}.$$ One can check that its cardinality
is $(d+1)^2.$ Consider the map
$\psi^d,$ the restriction of $\Psi^d$ to ${S_{m-1,n-1}}^{\Gamma_0}.$ 
\par
Denote by $MS^d$ the matrix of $\psi^d$ in the monomial bases. 
It is a
square matrix of size $(d+1)^2mn.$ 
\begin{remark}\
If $d=1,$ then $\psi^d=\phi$ and $MS^1=MP.$ For $d=2,$ 
 the matrices $MQ^2$ and $MS^2$ are 
denoted by $MQ$ and $MQ_w$ respectively in \cite{CGZ}.
\end{remark}
\begin{remark}\
It is straightforward to check that $MS^d$ is a maximal minor in 
$MQ^d.$  Furthermore, $\ker(MQ^d)$ is the $\KK$-vector space of 
moving
$d$-surfaces of bi-degree $(m-1,n-1)$ that follow the rational surface.
\end{remark}
\medskip
Consider the subset $\Gamma_1$ of $\Gamma_0$ defined by those $\gamma$ such that
$\gamma_1=0.$ Set
\begin{equation}
\label{tonta}
\rho^d:\, {S_{m-1,n-1}}^{\Gamma_1}\oplus S_{dm-1,dn-1} 
\to \, S_{(d+1)m-1,(d+1)n-1}
\end{equation}
the linear mapping which sends 
$\left((p_\gamma)_{\gamma\in\Gamma_1},q\right)$ to 
\begin{equation} 
\label{rho}
\sum_{\gamma\in\Gamma_1}p_\gamma x^\gamma + q\, x_1.
\end{equation}
Denote by $MT^d$ the matrix of $\rho^d$ in the monomial bases. 
One can check that it is a square matrix of the same size as $MS^d.$
\par
For a a square matrix $A,$ its determinant will be denoted by $|A|.$

\subsection{ Computing Resultants using Koszul Complexes}
We begin by reviewing the computation of  the determinant of short
exact sequence of  vector spaces 
(\cite{GKZ}, Appendix A). Consider the following exact complex
\begin{equation}
\label{scomplex}
0 \,\longrightarrow A \, 
\mathop{\longrightarrow}^{d_0} \, 
B 
\mathop{\longrightarrow}^{d_1}
\,C  \, \longrightarrow \, 0 \, . 
\end{equation}
Let $\{a_1,\ldots,a_p\}, \ \{b_1, \ldots, b_q\}, \ \{c_1,\ldots,c_r\}$
be  bases in $A, \, B,\, C$ respectively $(p+r=q).$
 In this case, 
the determinant of the complex with respect to these bases is equal to
the coefficient of proportionality
\begin{equation}
\label{wedge}
\frac{b_1\wedge\ldots\wedge b_q}{d_0(a_1)\wedge\ldots\wedge d_0(a_p)\wedge\widehat{c_1}
\wedge\ldots\wedge\widehat{c_r}}
\end{equation}
where $\widehat{c_1},\ldots \widehat{c_r}\in B$ satisfy $d_1(\widehat{c_i})
=c_i.$
\par
One can make an explicit computation of (\ref{wedge}) as follows:
\begin{enumerate}
\item \ Compute the matrices $D_0$ and $D_1$ corresponding to $d_0$ and $d_1$ 
respectively in the chosen bases.
\item \ Let $\overline{D_1}$ be the submatrix of $D_1$ given by all the $r$
rows and the first $r$ columns. Denote by $\overline{D_0}$ the submatrix of
$D_0$ given by the last $p$ rows and all the $p$ columns.
\item \ It turns out that $|\overline{D_0}|\neq0 \ \iff 
|\overline{D_1}|\neq0$ (cf. \cite{GKZ}). If this is the case, then
\begin{equation}
\label{mhdpl}
\det(\mbox{complex})=\frac{|\overline{D_1}|}{|\overline{D_0}|}
\end{equation}
\item If $|\overline{D_0}|=|\overline{D_1}|=0,$ then another maximal 
minor may be chosen in $D_0$ as follows: let 
$I=\{i_1,\ldots,i_r\}$ be an ordered subset of $r$ integers chosen from
$\{1,\ldots,q\}.$ Set
 $\overline{D_{0,I}},$ (resp. $\overline{D_{1,I}})$ the submatrix of
$D_0$ (resp. $D_1$) obtained by choosing all the columns (resp. rows) and the
rows (resp. columns) indexed by $\{1,\ldots,q\}\setminus I$ (resp. $I$). 
\par
 Change the order in the  basis of $B$ in order that the last $r$
elements are now indexed by $I.$ Using proposition $9$ in (\cite{GKZ},
Appendix 10) it is straightforward to check that
\begin{equation}
\label{eqcorr}
|\overline{D_{0,I}}|\,.\,\det(\mbox{complex})=(-1)^\sigma\,.\,|\overline{D_{1,I}}|,
\end{equation}
$\sigma$ being the parity of the permutation
$$\{i_1,\ldots,i_r,1,2,\ldots, i_1-1,i_1+1,\ldots, i_r-1,i_r+1,\ldots,q\}.$$
\end{enumerate}

Recall the definition of  $\Rd_{m,n}(f_1,f_2,f_3),$ 
the bihomogeneous resultant 
associated with a sequence of
three generic polynomials of bi-degree $(m,n)$~ 
(see for instance, \cite{dix,GKZ}): it is an irreducible polynomial in the 
coefficients of $f_i$ which vanishes after a specialization of the 
coefficients in a field $k$ 
if and only if the specialized system $f_i=0$ has a solution
in $\PP^1_\kk\times\PP^1_\kk.$
\par
One may compute  powers of
$\Rd_{m,n}(x_1,x_2,x_3)$ 
as a determinant of a three-term exact
Koszul complex as follows: consider the complex
of $\KK$-vector spaces:

\begin{equation}
\label{koszul}
\begin{array}{c}
0 \stackrel{}{\longrightarrow} { S_{m-1,n-1}}^2\oplus S_{(d-1)m-1,(d-1)n-1}  
\stackrel{\psi_0}{\longrightarrow}  
{S_{dm-1,dn-1}}^2\oplus S_{2m-1,2n-1}  
  \\
\stackrel{\psi_1}{\longrightarrow} S_{(d+1)m-1,(d+1)n-1}   \rightarrow  0  , 
\end{array}
\end{equation}
where $\psi_1$  and $\psi_0$
are the Koszul morphisms
\begin{equation}
\label{psi1}
\psi_1(p,q,r):= p\,x_1+q\,x_2+r\,x_3^{d-1};
\end{equation}
\begin{equation}
\label{psi0}
\psi_0\left(p,q,r\right):= \left(q\,x_3^{d-1}+r\,x_2, p\,x_3^{d-1}-r\,x_1,
-p\,x_2-q\,x_1\right).
\end{equation}
\medskip
\begin{proposition}
\label{tecn1}
The complex (\ref{koszul}) is  exact, and after a specialization
of the coefficients in a field $k$ it will remain exact  (as a complex of 
$k$-vector
spaces) if and
only if the bihomogeneous resultant of the specialized polynomials does not
vanish. The determinant of the complex
with respect to the monomial bases equals $\pm\Rd_{m,n}(x_1,x_2,x_3)^{d-1}.$ 
\end{proposition}
\begin{proof}\ 
Let $y_3(s,u;t,v)$ be a generic bihomogeneous polynomial of bi-degree 
$\left((d-1)m,(d-1)n\right).$ Consider the modified complex which is made
by replacing $x_3^{d-1}$ with $y_3$ in (\ref{psi1}), (\ref{psi0}).
\par
Because the polynomials $x_1,\ x_2$ and $y_3$ are bihomogeneous but do not have the
same bi-degree, the bihomogeneous resultant cannot be taken. However, 
there is another elimination operator  available: the {\it mixed resultant} 
associated with the sequence $(x_1,x_2,y_3)$ 
(\cite{GKZ}, Chapter $3$).
It is an irreducible polynomial in the 
coefficients of $x_1, x_2, y_3$ which vanishes if and only if these 
polynomials have a common root in $\PP^1\times\PP^1.$
\par
In order to compute it, we may apply the Cayley method for the study of resultants
(see \cite{GKZ}, Chapter $3$).  Let $\OO(d_1,d_2)$ denote the line bundle on 
$X:=\PP_{\kk}^1\times\PP_{\kk}^1$ whose sections are homogeneous
polynomials of degree $d_1$ in coordinates $(s:u)$ and degree $d_2$ in
coordinates $(t:v)$ on each $\PP^1.$
\par
Let $\LL_1=\LL_2=\OO(m,n)$ and $\LL_3=\OO((d-1)m,(d-1)n).$
Each $\LL_i$ is very ample, and we may regard polynomials of bi-degree 
$(m,n)$ with coefficients in $\kk$ as 
elements of $H^0(X,\LL_i),\, i=1,2,$ and polynomials of bidegree
$\left((d-1)m,(d-1)n\right)$ as belonging to $H^0(X,\LL_3).$
\par
Every specialization of $(x_1,x_2,y_3)$ defines a section $s$ of
the vector bundle
$E:=\LL_1\oplus\LL_2\oplus\LL_3.$  If we set $\MM:=\OO\left((d+1)m-1,(d+1)n-1\right),$
and construct the complex $C^{\bf\bullet}_{-}(\LL_1,\LL_2,\LL_3|\MM)$ 
(for a definition of this complex, see \cite{GKZ}, Chapter $3$), then we recover the modified complex (\ref{koszul}).
Moreover, it is not hard to check that this complex is stably twisted (i.e. has
no higher cohomology), so Proposition $4.1$ and Theorem $4.2$ of (\cite{GKZ}, Chapter $3$)
hold, and we have that the complex will be exact if and only if the resultant of the
specialized $(x_1,x_2,y_3)$ is not zero. Furthermore, the determinant of the
complex with respect to the monomial bases is equal to $\pm{\rm Res}(x_1,x_2,y_3).$
\par
The original complex  is recovered by specializing 
$y_3$ to $x_3^{d-1}.$  
Keeping in mind that 
$$ {\rm Res}(\tilde x_1,\tilde x_2, {\tilde x_3}^{d-1})=0\iff
\Rd_{m,n}(\tilde x_1,\tilde x_2, {\tilde x_3})=0,$$
we have that  the determinant of (\ref{koszul}) 
is equal to a power of $\Rd_{m,n}(x_1,x_2,x_3).$
\par
Comparing the degrees of both the bihomogeneous
and the mixed resultant in the coefficients of
$x_3$ and $y_3$ respectively (cf. \cite{CLO,GKZ}), we get that the determinant of the complex 
(\ref{koszul}) equals
$$ \pm \Rd(x_1,x_2,y_3)|_{y_3=x_3^{d-1}}= 
\pm\Rd_{m,n}(x_1,x_2,x_3)^{d-1}.$$
\end{proof}

Using the recipe given above, $\pm\Rd_{m,n}(x_1,x_2,x_3)^{d-1}$
may be computed using the 
following algorithm: 
\begin{enumerate}
\item \ Construct the matrices corresponding to the linear maps $\psi_0$ and 
$\psi_1$ with respect to the monomial bases;
\item \ Choose a  non vanishing maximal  minor $m_{1,I}$ in the matrix 
corresponding to 
$\psi_1,$  $I$ being a set of $(d+1)^2mn$ columns corresponding to vectors 
in the monomial 
basis of $\oplus_{i=1}^2 S_{dm-1,dn-1}\oplus S_{2m-1,2n-1}.$
\item \ Compute $m_{0,I},$ the maximal minor in the matrix representing 
$\psi_0,$ which consists of all rows {\bf not} indexed by $I.$
\item \ It turns out that $m_{0,I}\neq 0.$ Compute
$\frac{m_{1,I}}{m_{0,I}}.$ This quotient is equal to the determinant of the
complex.
\end{enumerate}
\medskip
\begin{remark}\
\label{last}
The following equality holds
for   every subset of indices $I:$
\begin{equation}
\label{pm}
m_{1,I} = \pm\Rd_{m,n}(x_1,x_2,x_3)^{d-1}\, m_{0,I}.
\end{equation}
\end{remark}

\subsection{The Relationship Between $|MS^d|$ and $|MP|$}
In order to prove the main result of this section, a preliminary lemma is
needed. Recall that $MT^d$ is the matrix associated with the linear map
(\ref{tonta}), and $MP$ is the matrix associated with (\ref{phi}).
\begin{lemma}
\label{auxlemma}
The following equality holds:
$$\pm |MT^d| =  |MP|^{d}\, \Rd_{m,n}(x_1,x_2,x_3)^{d(d-1)/2}.$$
\end{lemma}
\begin{proof}\
The proof will be  by induction on $d.$ For $d=1,$ it is clear than
$\rho^1$ and $\phi$ are the same functions, so the proposition follows
straightforwardly. 

Suppose then $d\geq2.$ The morphism $\rho^d$ may be factored 
 as follows:
{\small
\begin{equation}
\label{factrho}
{S_{m-1,n-1}}^{\Gamma_1}\oplus S_{dm-1,dn-1}
\mathop{\rightarrow}^{\psi_2}  {S_{dm-1,dn-1}}^2
\oplus S_{2m-1,2n-1}
\mathop{\rightarrow}^{\psi_1} S_{(d+1)m-1,(d+1)n-1}
\end{equation}}
where $\psi_1$ is the morphism defined in (\ref{koszul}) and
\begin{equation}
\label{psi2}
\psi_2\left(p_\gamma,q\right) =  \left(q,\,
\sum_{\gamma_2\geq1} p_\gamma x_2^
{\gamma_2-1}x_3^{\gamma_3}x_4^{\gamma_4}\,  , p_{(0,0,d,0)}x_3+p_{(0,0,d-1,1)}x_4 \right).
\end{equation}

Denote by $M_i$ the matrix corresponding to $\psi_i$ in the monomial bases
for $i=1,2.$
These are not square matrices (they have sizes
$(d+1)^2mn\times (2d^2+4)mn$ and $(2d^2+4)mn\times (d+1)^2mn$ respectively), 
but  applying the Cauchy-Binet  formula (see for instance \cite{HJ}), there is a relationship between their 
maximal minors and $|MT^d|:$
\begin{equation}
\label{bc} 
|MT^d| = \sum_{I} 
|M_{1,I}|\,|M_{2,I}|,
\end{equation}
the summation made over all sequences of integers $$I=(i_1,\ldots,i_{
(d+1)^2mn})$$ with $1\leq i_1\leq\ldots\leq i_{(d+1)^2mn}\leq (2d^2+4)mn,$ and
$M_{1,I}$ (resp. $M_{2,I}$) denotes the square submatrix of $M_1$
(resp. $M_2$) which is made by choosing the $(d+1)^2mn$ columns (resp. rows)
indexed by $I.$

\begin{remark}\
Note that $M_1$ is
the matrix corresponding to $\psi_1$ in the monomial bases  
and, for each $I,$ the maximal minor $m_{1,I}$ in step $2$ of the
algorithm outlined in the previous paragraph  is denoted 
$|M_{1,I}|$ in (\ref{bc}).
\end{remark}

\smallskip
Using remark \ref{last} and formula (\ref{bc}), one has
\begin{equation}
\label{sigue}
|MT^d| =\left( \Rd_{m,n}(x_1,x_2,x_3)\right)^{d-1} 
\sum_I (-1)^\sigma\,|M_{2,I}|\,m_{0,I}.
\end{equation}
An explicit computation of the $(2d^2+4)mn\times (d+1)^2mn$ matrix $M_2$ 
reveals the following structure:
\begin{equation}
\label{mp2}
M_2 =\left[\begin{array}{ccc}
\II & 0 & 0\\
0 &  B_1  & 0 \\
0 &  0  & B_2
\end{array}\right],
\end{equation} 
where $B_1$ and $B_2$ have sizes $d^2mn\times\frac{d(d+1)}{2}mn$ 
and $4mn\times 2mn$ respectively, and $\II$ denotes the identity matrix of size 
$d^2mn.$ 
\par
Gluing $M_2$ and the $(2d^2+4)mn\times ((d-1)^2+2) mn$ matrix $M_0$ 
corresponding to $\psi_0,$ 
one gets a square matrix $M:=\left[M_2,M_0\right]$ of size $2d^2+4$ 
which has the following structure:
\begin{equation}
\label{gmatrix}
M=\left[\begin{array}{ccccccc}
\II & 0 & 0&\vline&0&q\,x_3^{d-1}&r\,x_2\\
0 &  B_1  & 0 &\vline&p\,x_3^{d-1}&0&-r\,x_1\\
0 &  0  & B_2&\vline &-p\,x_2&-q\,x_1&0\\
\end{array}\right];
\end{equation} 
where the block $\left[p\,x_3^{d-1}\right]$ denotes the matrix 
corresponding to the 
linear map $S_{m-1,n-1}\to S_{dm-1,dn-1}$ which maps 
$p$ to $p\,x_3^{d-1},$ and the other blocks have the same meaning.
\par
It is easy to check that $\left[B_2,-p\,x_2, -q\,x_1\right]$ is a square
matrix. Moreover,
$$|\left[Q,-p\,x_2,- q\,x_1\right]| =\pm |MP|.$$
In the same way, it holds that 
$$|\left[B_1, -r\,x_1\right]|=\pm |MT^{d-1}|.$$
Then, the determinant of $M$ equals $\pm |MT^{d-1}||MP|,$ 
and  then the inductive hypothesis yields the following equality: 
\begin{equation}
\label{cremate}
|M|=\pm \left({\rm Res}_{m,n}(x_1,x_2,x_3)\right)^{(d-1)(d-2)/2}\,
|MP|^d.
\end{equation}
This determinant may also be computed as a sum
of maximal minors of  $M_2$ times their complementary minor in $M.$ This is 
exactly the sum which appears in (\ref{sigue}), i.e.
$$
|M| =\sum_I (-1)^\sigma\,|M_{2,I}|\,m_{0,I}.
$$
Replacing (\ref{cremate}) in (\ref{sigue}), the Lemma follows.
\end{proof}

\bigskip
\begin{theorem}
\label{mt}
$$\pm |MS^d| =  |MP|^{(d+1)d/2}\, \Rd_{m,n}(x_1,x_2,x_3)^{
(d+1)d(d-1)/6}.$$
\end{theorem}

\begin{proof}\
As in the Lemma, the proof will be  by induction on $d.$ 
For $d=1,$ it happens that
$\psi^1 =\phi,$ so the statement holds straightforwardly. 
\par
Take $d\geq2,$ and factor $\psi^d$ 
 as follows:
\begin{equation}
\label{factpsi}
 {S_{m-1,n-1}}^{\Gamma_0}
\mathop{\longrightarrow}^{\widetilde\psi_2} { S_{dm-1,dn-1}}^2
\oplus S_{2m-1,2n-1}
\mathop{\longrightarrow}^{\psi_1} S_{(d+1)m-1,(d+1)n-1}
\end{equation}
where 
\begin{equation}
\label{tildepsi2}
\widetilde\psi_2\left(p_\gamma\right) =  \left(\sum_{\gamma_1\geq1} 
p_\gamma x^\gamma,
\sum_{\gamma\in\Gamma_1,\gamma_2\geq1}p_\gamma x^\gamma  , 
p_{(0,0,d,0)}x_3+p_{(0,0,d-1,1)}x_4 \right).
\end{equation}
Denote by $\widetilde M_2$ the matrix corresponding to $\widetilde \psi_2$ 
in the monomial bases. 
The Cauchy-Binet formula gives the following relationship:
\begin{equation}
\label{bc2} 
|MS^d| = \sum_{I} 
|M_{1,I}|\,|\widetilde M_{2,I}|,
\end{equation}
where, as before, $I$ runs through all sequences of integers 
$(i_1,\ldots,i_{
(d+1)^2mn})$ satisfying $1\leq i_1\leq\ldots\leq i_{(d+1)^2mn}\leq (2d^2+4)mn,$ 
and $M_{1,I}$ (resp. $\widetilde M_{2,I}$) have the same meaning as in
(\ref{bc}).
\par
Proceeding as in the proof of the previous Lemma, one gets
\begin{equation}
\label{sigue2}
|MS^d| =\left( \Rd_{m,n}(x_1,x_2,x_3)\right)^{d-1} 
\sum_I (-1)^\sigma\,|\widetilde M_{2,I}|\,m_{0,I}.
\end{equation}
Gluing $\widetilde M_2$ and  matrix $M_0,$ 
one gets a square matrix
with the following structure:
\begin{equation}
\label{gmatrix2}
\widetilde M:=\left[\begin{array}{ccccccc}
MS^{d-1} & 0 & 0&\vline&0&q\,x_3^{d-1}&r\,x_2\\
0 &  B_1  & 0 &\vline&p\,x_3^{d-1}&0&-r\,x_1\\
0 &  0  & B_2&\vline &-p\,x_2&-q\,x_1&0\\
\end{array}\right].
\end{equation} 
Here, $B_1$ and $B_2$ are the same blocks which appear in (\ref{mp2}), 
and it is
easy to check that the block $\left[P,r\,x_1\right]$ is the matrix
$MT^{d-1}.$ 
 \par
Then, using inductive hypothesis and the previous Lemma, the following 
equalities hold: 
$$|\widetilde M|=\pm |MS^{d-1}|\,|MT^{d-1}|\,|MP|=
$$
$$
=\pm \left({\rm Res}_{m,n}(x_1,x_2,x_3)\right)^{(d-1)(d-2)(d+3)/6}\,
|MP|^{d(d+1)/2}.
$$
Computing this determinant as  sum
of maximal minors of  $\widetilde M_2$ times their complementary minor in $
\widetilde M,$ it appears the summation in 
(\ref{sigue2}). Replacing it  with this last expression, the Theorem follows
straightforwardly. 
\end{proof}
\bigskip
\begin{corollary}  {\rm  (Conjecture $(6.1)$ in \cite{CGZ})}
$$|MS^2| = |MP|^3\, {\rm Res}_{m,n} (x_1,x_2,x_3).$$
\end{corollary}
\smallskip
\begin{corollary}{\rm (general version of Theorem $4.1$ in \cite{CGZ})}
\par
Given  four bihomogeneous polynomials $\tilde x_1,\ \tilde x_2,\ \tilde x_3$ and $\tilde x_4$
of bi-degree $(m,n)$ and coefficients in $\CC.$ If ${\rm Res}_{m,n} (\tilde x_1,\tilde x_2,\tilde x_3)\neq0,$ then 
$|\tilde {MS^d}|=0$ implies
$|\tilde {MP}|=0.$  
\end{corollary}
\medskip
Here, $\tilde M$ means the matrix $M$ where the generic coefficients have been specialized with 
the coefficients of the $\tilde x_i.$
In the language of moving surfaces, Theorem \ref{mt} reads as follows:

\smallskip 
\noindent
{\it If $\Rd_{m,n}(\tilde x_1,\tilde x_2,\tilde x_3)\neq0,$ and there are
no moving planes of bi-degree\linebreak
$(m-1,n-1)$ which follow the rational surface, then the dimension of the $\CC$-vector
space of $d$-surfaces of bidegree $(m-1,n-1)$ that follow the rational surface 
is equal to $\frac{(d+1)d(d-1)}{6}\,m\,n.$}

\bigskip
\section{ The Validity of the Method of Moving Quadrics when the
Surface is not Properly parametrized}
\label{validity}
In this section,  we are going to discuss the validity of implicitization by moving quadrics with 
no base points, without requirements on the parametrization of the surface.
The main result will be an improvement of Theorem $4.2$ in \cite{CGZ}, extending the
validity of the method to the case when the parametrization is not generically one-to-one.
\par
We will set $d=2$ 
and $\KK=\QQ(c^i_{jk}).$
If $|MS^2|\neq0,$ then  $\ker (MQ^2)$ 
has dimension equal to $mn.$
Suppose without loss of generality that 
$MQ^2=\left[MS^2,R\right],$ 
where $R$ is a submatrix of $MQ^2$ of size $9mn\times mn.$
\par
In \cite{CGZ} (Proof of Theorem $4.2$),
a particular basis of the kernel of $MQ^2$ (i.e. a 
matrix $T$ of size $10mn\times mn$ such that $ MQ^2\cdot T=0$)
is considered:
$$ T:= \left[\begin{array}{c}
\bar T\\ \II
\end{array}
\right];$$ 
here, $\II$ denotes the identity matrix of order $mn.$
\par
Solving the linear equation $MQ\cdot T=0,$ one gets
$$T= \left[\begin{array}{c}
-{MS^2}^{-1}\cdot R\\ \II
\end{array}
\right].$$ 
$T$ has its rows indexed by the monomials $$s^i\,t^j x^\gamma,
\ 0\leq i\leq m-1, \ 0\leq j\leq n-1, |\gamma|=2,$$ and the last $mn$ rows
corresponds to the monomials $$s^i\,t^j x_4^2,
\ 0\leq i\leq m-1, \ 0\leq j\leq n-1.$$
\par
Consider $\TT,$ the matrix which results reordering the rows of $T$ as
follows:
\begin{equation}
\label{orden}
\begin{array}{c}
x_1^2 ,\ x_2^2,\ x_3^2,\ x_4^2,\ x_1\,x_2,\ x_1\,x_3,\ x_1\,x_4,\ x_2\,x_3
,\ x_2\,x_4 ,\ x_3\,x_4, \\
s\left(x_1^2 ,\ x_2^2,\ x_3^2,\ x_4^2,\ x_1\,x_2,\ x_1\,x_3,\ x_1\,x_4,\ x_2\,
x_3,\ x_2\,x_4, \ x_3\,x_4\right), \\ 
s\,t\left(x_1^2 ,\ x_2^2,\ x_3^2,\ x_4^2,\ x_1\,x_2,\ x_1\,x_3,\ x_1\,x_4,\ 
x_2\,x_3,\
x_2\,x_4 ,\ x_3\,x_4\right),\\
s\,t^2\left(x_1^2 ,\ x_2^2,\ x_3^2,\ x_4^2,\ x_1\,x_2,\ x_1\,x_3,\ x_1\,x_4,\ 
x_2\,x_3,\
x_2\,x_4,\ x_3\,x_4\right),\\
s\,t^3\left(x_1^2 ,\ x_2^2,\ x_3^2,\ x_4^2,\ x_1\,x_2,\ x_1\,x_3,\ x_1\,x_4,\ 
x_2\,x_3,\
x_2\,x_4,\ x_3\,x_4\right),\\
\ldots \\
\end{array}
\end{equation}
The third row of $\TT,$ for example, corresponds to the 
monomial $x_3^2;$ its eleventh row is indexed by $s\,x_1^2.$
\par
Let $X_1,\ X_2,\ X_3$ and $X_4$ be indeterminates over $\KK.$ Consider the
following vector in $\ZZ\left[X_1,X_2,X_3,X_4\right]^{10}:$
\begin{equation}
\label{c} 
C:=\left(X_1^2, \, X_2^2,\, X_3^2,\, X_4^2,\, X_1\, X_2,\, X_1\,X_3,\, 
X_1\,X_4,\, X_2\,X_3,
\, X_2\,X_4, \, X_3\,X_4 \right).
\end{equation}
Let $M\in\ZZ[X_1,X_2,X_3,X_4]^{mn\times 10mn}$ be the matrix defined 
as follows:
$$ M:= \left[
\begin{array}{ccccc}
C&C&C&\ldots&C\\
0&C&C&\ldots&C\\
0&0&C&\ldots&C\\
&&\ldots&&\\
0&0&0&\ldots&C
\end{array}
\right].$$
Recall the construction given  in the proof of Theorem $4.2$ in
\cite{CGZ}
for computing the implicit equation of the parametric surface:
each column $T_{\alpha,\beta}$ of $\TT$ encodes the coordinates (
in the monomial basis) of a 
moving quadric of bidegree
$(m,n)$ which follows the surface. Write $T_{\alpha,\beta}$ as follows:
$$T_{\alpha,\beta} = \sum_{ i \leq m-1,\ j\leq n-1 } T^{\alpha,\beta}_{i,j}
\, s^i\, t^j,$$
where 
$$T^{\alpha,\beta}_{i,j}= \sum_{|\gamma|=2} a^\gamma_{i,j} x^\gamma.$$
Set $\widetilde T^{\alpha,\beta}_{i,j}:= 
\sum_{|\gamma|=2} a^\gamma_{i,j} X^\gamma,$
and consider  the square matrix 
$$\widetilde T:= \left[T^{\alpha,\beta}_{i,j}\right]
\in \KK[X_1,X_2,X_3,X_4]^{mn\times mn},$$
where its rows are indexed by $(i,j)$ and columns by $(\alpha,\beta).$

\begin{remark}\
$\widetilde T$ is the matrix called $M$ in the proof of Theorem $4.2$ in
\cite{CGZ}.
\end{remark}

\begin{proposition}
\label{formal} \
$|\,\widetilde{T}\,| =  \pm|M\cdot\TT|. $
\end{proposition}
\begin{proof}\
Set $\R:= M\cdot\TT.$
Computing explicitly the last row of 
$M\cdot\TT,$  because of the order given to the rows of $\TT$ in (\ref{orden}),
one obtains:
$$\left[\widetilde T_{m-1,n-1}^{\alpha,\beta}\right]_{\alpha,\beta},$$ 
which coincides with one of the rows in $\widetilde T.$
\par
In the same way, one can check that the row inmediately before the last
in $\R,$ is the following
$$ \left[\widetilde T_{m-1,n-2}^{\alpha,\beta}+
\widetilde T_{m-1,n-1}^{\alpha,\beta}\right]_{\alpha,\beta},$$  
so, substracting from it the last row, one gets another row of $\widetilde T.$
\par
A similar situation happens in all rows. This implies that the
matrix $\R$ may be transformed in the matrix $\widetilde T$ applying
operations on its rows which do not change the determinant.  
\end{proof}

\medskip
The following Proposition will be useful in the sequel.
\begin{proposition}
\label{4.3}
Set $$P(X_1,X_2,X_3):=\Rd_{m,n}\left(x_1-X_1\,x_4,
x_2-X_2\,x_4,x_3-X_3\,x_4\right).$$
Then, $P(X)$ is an irreducible polynomial in $\KK[X_1,X_2,X_3]$
and its degree in the variables $X_i$ is equal to $mn.$
\end{proposition}
\begin{proof}\
The fact that the degree of $P$ in $X_i$ is $mn$ can be easily checked in
Dixon's matrix (cf. \cite{dix}).
\par
Set ${\bf Z}:= \ZZ[c^i_{jk}, X_1,X_2,X_3].$
Suppose 
$$P= A(c,X)\,.\,B(c,X)$$
where $A,B\in {\bf Z}.$ The polynomial $P$ is homogeneous in the variables $c^i_{ij},$
which implies that $A$ and $B$ are also homogenous in the $c^i_{jk}.$
\par
Specializing $X_i\mapsto0,$ one has that
$$\Rd_{m,n}(x_1,x_2,x_3)= A(c,0)\,.\,B(c,0).$$
But the left-hand side is irreducible. This implies that one of the 
factors must have degree $0$ in the variables $c^i_{jk},$ lets say $B.$
The factorization now reads as follows
$$P= A(c,X)\,.\,B(X),$$
where $B\in\ZZ[X_1,X_2,X_3].$ But if $\deg_X{B}\geq1,$ then
the variety \linebreak 
$\left(B=0\right)\neq\emptyset$ in $\CC^3.$
On the other hand, it is well known that, for a given 
family of bihomogeneous polynomials $\tilde x_i(s,u;t,v)$ with 
no base points,
the equation $$\Rd_{m,n}\left(\tilde x_1-X_1\,\tilde x_4,
\tilde x_2-X_2\,\tilde x_4,\tilde x_3-X_3\,\tilde x_4\right)=0,$$
is a power of the implicit equation of the rational surface defined by
$\left(\frac{\tilde x_1}{\tilde x_4},\frac{\tilde x_2}{\tilde x_4},
\frac{\tilde x_3}{\tilde x_4}\right)$
(actually, if we use Dixon matrices for computing this resultant, we will find moving planes
encoded in their rows; see the comment at the end of Example $3$ in \cite{CGZ}).
If $B$ had a nonempty zero locus, this would imply that  {\it every}
rational parametric surface of bi-degree $(m,n)$ without
base points will contain the zero locus
of $B,$ which is impossible.
\end{proof}
\begin{corollary}
\label{4.4}
Let $P^{h}(X_1,X_2,X_3,X_4)$ be the homogeneization of $P$ up to degree
$mn.$ Then, $P^h$ is irreducible in $\ZZ[c^i_{jk},X_1,X_2,X_3,X_4].$
\end{corollary}
The following may be regarded as the main result of this section:
\begin{theorem}
\label{mth}
$$ \Rd_{m,n}(x_1,x_2,x_3)\,.\,|M\cdot\TT| =  P^h,$$
\end{theorem}
\begin{proof}\
To begin, it will be proven that the ratio between
$|M\cdot\TT|$ and $P^h$ is in $\KK.$ 
Note that there is a positive integer $z$ such that
$$\R^{'}:=|MS|^z\,.\,\det(M\cdot\TT)\in
\ZZ[c^i_{jk},X_1,X_2,X_3,X_4]$$
(this is due to the fact that, in the construction of $\TT,$ one needs 
the inverse of $MS.$)
\par
Set $\R^{''}:=|MS|\,.\,\R^{'}.$
Both $R^{''}$ and $P^h$ are homogeneous in $X$ and of the same degree, 
which is equal to
$mn,$ so it will be enough to show that one of them is a
polynomial multiple of the other.
\par
Suppose that the variables have been specialized:
$$\left(c^i_{jk},X\right)\mapsto\left(\tilde c^i_{jk},\widetilde X\right)$$ 
with the following conditions:
\begin{enumerate}
\item \ $|\widetilde MS|\neq0,$ where $\widetilde MS$ denotes the matrix
$MS$ after specializing $c\mapsto \tilde c.$
\item \ $P^h(\tilde c^i_{jk},\widetilde X)$ vanishes. 
\end{enumerate}
This means that the projective point $\left(\widetilde X_1:\widetilde X_2:
\widetilde X_3:\widetilde X_4\right)$ belongs to the rational surface defined by
$\left(\frac{\tilde x_1}{\tilde x_4},\frac{\tilde x_2}{\tilde x_4},
\frac{\tilde x_3}{\tilde x_4}\right).$ 
\par
Using  $\widetilde MS\neq 0$ and the same
argument as in the proof of Theorem $4.2$ in \cite{CGZ}, 
one can verify that, if $\left(\widetilde X_1:\widetilde X_2:
\widetilde X_3:\widetilde X_4\right)$ belongs to that rational surface, 
the determinant of $\widetilde T$ must vanish. Proposition 
\ref{formal} implies that $\R^{''}$ must vanish.
Hence,
$$ \left(P^h=0\right)\subset \left(\R^{''}=0\right) \ \ \ \mbox{if}\ 
|MS|\neq0.$$
If $|MS|=0,$ the inclusion holds trivially and, using
the Hilbert's Nullstellensatz, one can conclude that a power of 
$\R^{''}$ must be a multiple of $P^h.$
As $P^h$ is irreducible,  $\R^{''}$ must be a multiple
of $P^h,$ so
$$|M\cdot\TT|= c\,.\, P^h, $$
where $c\in\KK$ as expected. 
\par
In order to compute $c,$ do the following replacement:
$$X_1,\ X_2,\ X_3\mapsto 0, \ \ X_4\mapsto 1.$$ 
Then, it will hold that
$$\begin{array}{ccc}
P^h&\mapsto &\Rd_{m,n}(x_1,x_2,x_3)\\
 \widetilde T &\mapsto &\II_{mn}
\end{array}$$
if the columns of $\widetilde T$ are properly ordered 
(cf. \cite{CGZ}, proof of Theorem $4.2$).
So, $c$ must be $\pm\frac{1}{\Rd_{m,n}(x_1,x_2,x_3)}$
\end{proof}

\begin{corollary}
Given a family of bihomogeneous polynomials 
\linebreak
$\tilde x_i(s,u;t,v)\in\CC[s,u,t,v], 
\ 1\leq i \leq 4.$
Suppose that the parametrization $$\left(\frac{\tilde x_1}{\tilde x_4},
\frac{\tilde x_2}{\tilde x_4},\frac{\tilde x_3}{\tilde x_4}\right)$$ 
defines a surface which has no base points, and that
there are no moving planes of bi-degree $(m-1,n-1)$ following the
surface. Then, the method
of moving quadrics always computes a power of the implicit equation of 
this surface.
\end{corollary}
\begin{proof}\
In order to make the method work correctly,
one can suppose without loss of generality that $\Rd_{m,n}(\tilde x_1,\tilde x_2,\tilde x_3)
\neq0.$

The method computes $|\widetilde T| = \pm |M\cdot\TT|$ which is 
equal to a constant times $\tilde {P^h},$ because of
Theorem
\ref{mt}. But $\tilde {P^h}$ vanishes if and only if the projective point 
$(X_1,X_2,X_3,X_4)$ belongs to the implicit surface.

This, combined with the fact that the implicit equation of the parametric 
surface is always irreducible, completes the proof.
\end{proof}

\bigskip
\section{The triangular surface case}
\label{tres}
In the triangular case, similar results hold. 
We shall denote polynomials, linear maps and
its matrices as in Section \ref{dos}, though the reader should keep in
mind that everything will now be homogeneous rather than bi-homogeneous.
\par
Let
$$ x_i\left(s,t,u\right) = \sum_{j+k\leq n}
c^i_{jk} \ s^j\,t^k\,u^{n-j-k} \ \ \ \ i=1,\ldots,4$$
be four generic polynomials in three variables of degree $n.$ 
\par
Set, as before, $\KK:= \QQ(c^i_{jk}),$
and let $S_{l}$ be the space of homogeneous polynomials in three variables
of degree  $l,$ with coefficients in $\KK.$
\par
Consider now:
\begin{equation}
\begin{array}{cccc}
\phi:& {S_{n-1}}^4& 
\to & S_{2n-1}\\
&\left(p_1,p_2,p_3,p_4\right)&\mapsto & \sum_{i=1}^4 p_i\, x_i,\\
\end{array}
\end{equation}
and, let $MP$ be the matrix of $\phi$ in the
monomial bases (cf. \cite{CGZ}). 
\par
$MP$ has size $(2n^2+n) \times (2n^2+2n).$ In order to have a square submatrix,
let $\ii\subset\{(j,k): \, j+k\leq n\}$  with $|\ii|=n.$ Define
$MP_{\ii}$ by removing the $n$ columns in $MP$
corresponding to $$s^i\,t^j\,u^{n-i-j}\,x_4, \ (i,j)\in\ii.$$
As in the  Section \ref{dos}, $MP_{\mathcal I}$ corresponds to a maximal square submatrix
of the coefficient matrix of the system generated by the moving planes
of degree $n-1$ that follow the rational surface.
\par
Set  $\Gamma$ and $\Gamma_0$ as before, 
and consider the maps
$$\Psi^d:\, {S_{n-1}}^{\Gamma}\, 
\to \, S_{(d+1)n-1}$$ 
which sends
the sequence $\left(p_\gamma\right)_{\gamma\in\Gamma}$ to 
$
\sum_{\gamma\in\Gamma}p_\gamma x^\gamma,
$
and
$$\psi^d:\, {S_{n-1}}^{\Gamma_0}\, 
\to \, S_{(d+1)n-1},$$ 
its restriction.
Let $MQ^d$  and $MS^d$ be the matrices of $\Psi^d$ and 
$\psi^d$ in the monomial bases, respectively.
$MS^d$ has size  $\frac{((d+1)n+1)(d+1)n}{2}\times
\frac{(d+1)^2n(n+1)}{2}.$ 
To get a square submatrix of it,  remove all the columns corresponding
to 
$$
s^i\,t^j\,u^{n-i-j}\,x^\gamma, \ (i,j)\in \ii,\ \gamma_4=1,
$$
and define $MS^d_{\ii}$ to be the remaining submatrix. It is associated to
a linear map
\begin{equation}
\label{psidi}
\psi_{\ii}^d:\, {S_{n-1,\ii}^{\Gamma_0\setminus\Gamma_4}}\oplus{S_{n-1}}^{\Gamma_4}\, 
\to \, S_{(d+1)n-1}
\end{equation}
where $S_{n-1,\ii}$ is the $\KK$-vector space generated by all
monomials whose exponents are not in $\ii,$ and $\Gamma_4$ is the set of all 
multiindices $\gamma\in\Gamma_0$ such that $\gamma_4=0.$
Also, $MS^d_{\ii}$ may be regarded as  a maximal square submatrix of the
coefficient matrix of the moving $d$-surfaces of degree $n-1$ that follow the rational surface.
\par
Let $\Rd_{n}(f_1,f_2,f_3),$  be
the multivariate resultant 
associated with a sequence of
three generic polynomials of degree $n$~ 
(cf. \cite{dix,CLO,GKZ}): it is an irreducible polynomial in the 
coefficients of $f_i,$ which vanishes after a specialization of the 
coefficients in an algebraically closed field $\kk$ 
if and only if the specialized system $f_i=0$ has a solution
in $\PP^2(\kk).$
\par
In this situation, a similar result holds:
\begin{theorem}
\label{mtp}
$$|MS^d_{\ii}|= \pm |MP_{\ii}|^{(d+1)d(d-1)/6}\, 
\left(\Rd_n(x_1,x_2,x_3)\right)^{(d+1)d/2}.$$
\end{theorem}

\begin{proof}\
The proof follows applying mutatis mutandis all the tools developed in 
Section \ref{dos}.
Consider the
following Koszul Complex: 
\begin{equation}
\label{koszul2}
0 \,\longrightarrow  {S_{n-1}}^2\oplus{S_{(d-1)n-1}} \, 
\mathop{\longrightarrow}^{\psi_0} \, 
{S_{dn-1}}^2\oplus S_{2n-1}  
\mathop{\longrightarrow}^{\psi_1}
\, S_{(d+1)n-1}  \, \longrightarrow \, 0 \, .
\end{equation}
Here, $\psi_1$ and $\psi_0$ are defined by (\ref{psi1}) and (\ref{psi0})
respectively. A similar version of Proposition \ref{tecn1} holds, applying
the same trick used there to the formulation of the multivariate resultant
as the determinant of a Koszul Complex (cf. \cite{cha,de}).
\begin{proposition}
\label{tecn2}
The complex (\ref{koszul2}) is  exact, and after a specialization
of the coefficients $\tilde{x}_1,\, \tilde{x}_2,\,\tilde{x}_3,$
it will be exact if and
only if  $\Rd_{n}(\tilde{x}_1,\tilde{x}_2,\tilde{x}_3)$ does not vanish.
The determinant of the complex
with respect to the monomial bases is equal to $\pm\Rd_{n}(x_1,x_2,x_3)^{d-1}.$ 
\end{proposition}

The linear map $\psi_{\ii}^d$ may be factored as follows:
\begin{equation}
 {S_{n-1,\ii}^
{\Gamma_0\setminus\Gamma_4}}\oplus{S_{n-1}}^{\Gamma_4}
\mathop{\longrightarrow}^{\widetilde\psi_2} \, 
{S_{dn-1}}^{2}\oplus S_{2n-1}  
\mathop{\longrightarrow}^{\psi_1}
\, S_{(d+1)n-1},  
\end{equation}
where $\widetilde \psi_2$ is defined as in (\ref{tildepsi2}). 
\par
In order to apply to this situation the proof of Theorem (\ref{mt}), 
the ``triangular'' version of Lemma (\ref{auxlemma}) is 
needed. Let $$\rho^d:\, {S_{n-1}}^{\Gamma_1}\oplus S_{dn-1} 
\to \, S_{(d+1)n-1}$$ 
be the linear map defined as in (\ref{rho}), and set $MT^d_{\ii}$ the
matrix in the monomial bases corresponding to the restriction of $\rho^d$
to $$\left(
{S_{n-1,\ii}^{\Gamma_1\setminus\Gamma_4}}\oplus{S_{n-1}}^{\Gamma_1\cap
\Gamma_4}\right)\oplus S_{dn-1}.$$ It is a square matrix of the same size
as $MS^d_{\ii}.$ The following equality follows straightforwardly:

\begin{lemma}
$$\pm |MT^d_{\ii}| =  |MP_{\ii}|^{d}\, \Rd_{n}(x_1,x_2,x_3)^{d(d-1)/2}.$$
\end{lemma}
Now, the same proof performed in Theorem (\ref{mt}) for the tensor product
case, may be applied for triangular polynomials, giving the desired result.
\end{proof}

\bigskip
\begin{corollary}  {\rm  (Conjecture $(6.2)$ in \cite{CGZ})}
\par
$$|MS^2_{\ii}| = |MP_{\ii}|^3\, {\rm Res}_{n} (x_1,x_2,x_3).$$
\end{corollary}
\smallskip
\begin{corollary}{\rm (general version of Theorem $5.1$ in \cite{CGZ})}
\par
If ${\rm Res}_{n} (x_1,x_2,x_3)\neq0,$ then $|MS_{\ii}^d|=0$ implies
$|MP_{\ii}|=0.$
\end{corollary}

Theorem \ref{mtp} has the following interpretation in the language of moving surfaces:

\smallskip
\noindent
{\it If $\Rd_{n}(\tilde x_1,\tilde x_2,\tilde x_3)\neq0,$ and there are
exactly  $n$ linearly independent moving planes of degree
$n-1$ that follow the rational surface, then the dimension of the $\CC$-vector
space of $d$-surfaces of degree $n-1$ that follow the rational surface is  
equal to $\frac{n(d+1)d(dn+d+5-n)}{12}.$}

\begin{proof}\
The fact that there are $n$ linearly independent moving planes of degree $n-1,$ implies that $MP$ has
maximal rank. Lemma $5.2$ in \cite{CGZ} implies that there exists in index set $\ii, \ |\ii|=n,$ such that
$MP_{\ii}$ is non-singular.
\end{proof}

\subsection{ The Method of Moving Quadrics when the Surface is not Properly 
parametrized}
We are going to give here the analogue of Section \ref{validity} for triangular surfaces.
The main result will be an improvement of Theorem $5.2$ in \cite{CGZ}, extending the
validity of the method to the case when the parametrization is not generically one-to-one.
Some extra care must be taken, because the method combines moving planes and moving 
quadrics . 
\par As before, set $d=2,$ and 
suppose that $\ii$ is fixed. If $|MS_{\ii}^2|\neq0,$ the kernel of 
$MQ^2,$ i.e. the space of moving quadrics which follow the surface, 
will have dimension equal to $\frac{n^2+7n}{2}.$
\par
As in Section \ref{validity}, suppose that
$MQ^2=\left[MS_{\ii}^2,R_{\ii}\right],$
where $R_{\ii}$ has $\frac{n(n+1)}{2}+3n$ columns.
In \cite{CGZ} (proof of Theorem $5.2$),
$(n^2-n)/2$ linearly independent vectors are chosen from the basis of 
 the kernel of $MQ^2$ by considering  a
matrix $T$ of 
$(5n+5n^2)\times (n^2-n)/2$ such that $ MQ^2\cdot T=0$ 
as follows:
$$ T:= \left[\begin{array}{c}
\bar T\\ 
0\\
\II
\end{array}
\right].$$ 
Here, $0$ is a block of $4n$ rows indexed by 
 $s^i\,t^j\,x_i^a x_4^b$, $\  (i,j)\in \ii,\, b\geq1,$ and
$\II$ denotes the identity matrix of order $(n^2-n)/2$ whose rows are indexed
by $s^i\,t^j\,x_4^2 \ (i,j)\notin \ii.$
\par
With this structure, one can check that 
$$T= \left[\begin{array}{c}
-{MS_{\ii}^2}^{-1}\cdot \widetilde{R_{\ii}}\\
0\\ 
\II
\end{array}
\right],$$
where $\widetilde{R_{\ii}}$ denotes a block of $R_{\ii}.$
\par
Because of Theorem \ref{mtp},
the fact that $|MS^2_{\ii}|\neq0$ implies that $MP^2_{\ii}$ is not
singular, hence one can also find $n-1$ linearly independent
moving planes that follow the
surface by solving the system $MP\cdot T'=0,$ where
$MP = \left[MP_{\ii}\,,\,R'\right]$ and 
$$T'= \left[\begin{array}{c}
-{MP_{\ii}}^{-1}\cdot \widetilde{R'}\\ 
\II
\end{array}
\right].$$ 
As in the previous section, let $\TT$  be the matrix made from
$T$ by ordering its rows with the order defined in (\ref{orden}). 
\par
Consider also the following order:
$$\begin{array}{c}
x_1 ,\ x_2,\ x_3,\ x_4, \\
s\left(x_1 ,\ x_2,\ x_3,\ x_4\right), \\ 
s\,t\left(x_1 ,\ x_2,\ x_3,\ x_4\right),\\
s\,t^2\left(x_1 ,\ x_2,\ x_3,\ x_4\right),\\
\ldots \\
\end{array}
$$
and let $\TT'$ be the matrix made from $T'$ by ordering its
rows with it.
\par
Let  $C\in \ZZ\left[X_1,X_2,X_3,X_4\right]^{10}$
be the vector defined in (\ref{c}), and set $C':=(X_1,X_2,X_3,X_4).$
Consider also $M\in\ZZ[X_1,X_2,X_3,X_4]^{n(n+1)/2\times 10n(n+1)/2}$ 
and $M'\in\ZZ[X_1,X_2,X_3,X_4]^{n(n+1)/2\times 4n(n+1)/2}$ 
as follows:
$$ M:= \left[
\begin{array}{ccccc}
C&C&C&\ldots&C\\
0&C&C&\ldots&C\\
0&0&C&\ldots&C\\
&&&\ddots&\\
0&0&0&\ldots&C
\end{array}
\right],$$
$$ M':= \left[
\begin{array}{ccccc}
C'&C'&\ldots&C'\\
0&C'&\ldots&C'\\
&\ldots&&\\
0&0&\ldots&C'
\end{array}
\right].$$
\par
Let $\widetilde T$ be the square matrix of size $n(n+1)/2$
constructed in the proof of Theorem $5.2$ in
\cite{CGZ} by collecting all the coefficients of these moving planes and
moving quadrics. The following Proposition may be proven mutatis
mutandis the arguments given in Proposition \ref{formal}:

\begin{proposition}
$$|\widetilde T|= \pm |\left[M\cdot\TT\, , 
\, M'\cdot\TT'\right]|. $$
\end{proposition}
\medskip
Similarly, one can formulate and prove versions of Proposition \ref{4.3}, Corollary \ref{4.4} and
Theorem \ref{mth} for the triangular case. This leads to the following corollary:
\begin{corollary}
Suppose that the surface $$\left(\frac{\tilde x_1}{\tilde x_4},
\frac{\tilde x_2}{\tilde x_4},\frac{\tilde x_3}{\tilde x_4}\right)$$ 
has no base points, and that
there are exactly $n$  moving planes of degree $n-1$ following the
surface. Then, the method
of moving quadrics always compute a power of the implicit equation of 
the surface.
\end{corollary}

\bigskip
\noindent {\bf Acknowledgements:} 
The author would like to express his deep gratitude to his advisor 
Alicia Dickenstein,
not only for  her
thorough reading of preliminary drafts of this paper and very
thoughtful suggestions for improvement, 
but for introducing him to this interesting area. 
He is also grateful
to David Cox and the referees for their very helpful suggestions.

\bigskip

\bigskip
\noindent\email{cdandrea@dm.uba.ar}

\end{document}